\documentclass[12pt,oneside]{amsart}
\usepackage[english]{babel}
\usepackage{amssymb}
\usepackage{amsthm}
\usepackage[foot]{amsaddr}
\usepackage{amsmath}
\usepackage{a4wide}
\usepackage{esvect}
\usepackage{hyperref}
\usepackage{verbatim}
\usepackage[dvips]{graphicx}
\usepackage{t1enc}
\usepackage{color}
\usepackage{geometry}

\usepackage{tikz}
\usepackage{caption}
\usepackage{subcaption}
\newtheorem{thm}{Theorem}

\newtheorem{conjecture}[thm]{Conjecture}

\newtheorem{corollary}[thm]{Corollary}

\newcommand{\ignore}[1]{}

\begin{document}

\author{Ervin Gy\H ori}
\address[Ervin Gy\H ori]{Alfr\'ed R\'enyi Institute of Mathematics / Central European University}
\email[Ervin Gy\H ori]{gyori.ervin@renyi.mta.hu}
\thanks{Research  of the first author was supported by OTKA grant 116769.}

\author{Tam\'as R\'obert Mezei}
\address[Tam\'as R\'obert Mezei]{Central European University}
\email[Tam\'as R\'obert Mezei]{tamasrobert.mezei@gmail.com}

\author{G\'abor M\'esz\'aros}
\address[G\'abor M\'esz\'aros]{Alfr\'ed R\'enyi Institute of Mathematics}
\email[G\'abor M\'esz\'aros]{meszaros.gabor@renyi.mta.hu}
\thanks{Research of the third author was supported by OTKA grant 116769.}

\title[Note on Terminal-Pairability in Complete Grid Graphs] 
{Note on Terminal-Pairability in Complete Grid Graphs}

\linespread{1.3}
\pagestyle{plain}

\begin{abstract}
 We affirmatively answer and generalize the question of Kubicka, Kubicki and Lehel \cite{grid} concerning the path-pairability of high-dimensional complete grid graphs. As an intriguing by-product of our result we significantly improve the estimate of the necessary maximum degree in path-pairable graphs, a question originally raised and studied by Faudree, Gy\'arf\'as, and Lehel \cite{pp}.
\end{abstract}
\maketitle 

\section{Introduction}
We discuss a graph theoretic concept of {\it terminal-pairability} emerging from a practical networking problem introduced by Csaba, Faudree, Gy\'arf\'as, Lehel, and Shelp \cite{CS} and further studied by Faudree, Gy\'arf\'as, and Lehel \cite{mpp,F,pp} and by Kubicka, Kubicki and Lehel \cite{grid}. A graph $G$ on $2k$ vertices is called {\it path-pairable} if,
for any pairing of the vertices
$X = \{x_1,\dots,x_k\}$ and $Y = \{y_1,\dots,y_k\}$ of $G$, there exist $k$
edge-disjoint $x_iy_i$-paths joining the pairs. The vertices of the set $X\cup Y$ are often called {\it terminals} while the pairs
$(x_i,y_i)$ of terminals are simply called {\it pairs}.

A more general pairability concept including the above discussed one is {\it terminal-pairability.} Let $G$ be a graph with maximum degree $\Delta$ and with vertex set $V(G) = T(G)\cup I(G)$ where the set $T(G)$ consists of $t$ ($t$ even) vertices of degree 1. We call $G$ a {\it terminal-pairable} network if for any pairing of the vertices of $T(G)$ there exist edge-disjoint paths in $G$ between the paired vertices. $T(G)$ is referred to as the set of {\it terminal nodes} or {\it terminals} and $I(G)$ is called the set of interior nodes of the network.  For an inner vertex $v$ we denote the number of terminal vertices incident to $v$ by $d_{T(G)}(v)$. Observe that for a graph $G$ being path-pairable is equivalent to being terminal-pairable with $d_{T(G)}(v)= 1$ for every $v\in V(G)$.

Given a terminal-pairable network $G$ with a particular pairing $\mathcal{P}$ of the terminals the {\it demand (multi)graph} $D=(V(D),E(D))$ is defined as follows: we set $V(D)=I(G)$ and join two vertices of $V(D)$ by an edge if the corresponding vertices in $I(G)$ are incident to terminals that form a pair in $\mathcal{P}$. Obviously, $E(D)=\frac{|T(G)|}{2}$ and $d(v)= d_{I(T)}(v)$ for every $v\in V(D)$, thus in fact $\Delta(D)= \max_{v\in I(G)} d_{I(T)}(v)$. 
Observe also that a terminal pairing problem is fully described by the underlying network $G$ and the demand graph $D$.

A long-standing open question concerning path-pairability of graphs is the
minimal possible value of the maximum degree $\Delta(G)$ of a path-pairable
graph $G$. Faudree, Gy\'arf\'as, and Lehel \cite{pp} proved that the
maximum degree has to grow together with the number of vertices in path-pairable
graphs. They in fact showed that a path-pairable graph with maximum degree
$\Delta$ has at most $2\Delta^\Delta$ vertices. The result yields a lower
bound of order of magnitude $\frac{\log n}{\log \log n}$ on the maximum degree of
a path-pairable graph on $n$ vertices. This bound is conjectured to be
asymptotically sharp, although to date only constructions of much higher order of magnitude have been found. The best
known construction is due to Kubicka, Kubicki, and Lehel
\cite{grid} who showed that two dimensional complete grids on an even number of vertices are path-pairable. A two dimensional complete grid is the {\it Cartesian product} $K_s\square K_t$ of two complete graphs $K_s$ and $K_t$ and it can be constructed by taking the Cartesian product of the sets $\{1,2,\dots\,s\}$ $\{1,2,\dots\,t\}$ and joining two vertices if they share a coordinate. Higher dimensional complete grids can be defined similarly; for a more detailed introduction of the Cartesian product of graphs we refer the reader to \cite{product}.

With $s=t$ the construction of Kubicka, Kubicki, and Lehel gives examples of path-pairable graphs on $n=s\cdot t$ vertices with maximum degree $2\cdot\sqrt{n}$. This bound was recently improved to $\sqrt{n}$ by M\'esz\'aros \cite{me_pp}. It was also conjectured in \cite{grid} that $K_t\square K_t\square K_t$ is path-pairable for sufficiently large even values of $t$.

In this paper we significantly improve the upper bound on the minimal value of $\Delta$ by proving path-pairability of high dimensional complete grids. We eventually study the more general terminal-pairability variant of the above path-pairability problem and prove the following Theorem:

\begin{thm}\label{main}
Let $G=K_t^n $ and let $D=(V(D),E(D)$ be a demand graph with $V(D)=V(K_t^n)$ and $\Delta(D)\leq \lfloor\frac{t}{6}\rfloor- 1$ even. Then every demand edge of $D$ can be assigned a path in $G$ joining the same endpoints such that the system of paths is edge-disjoint. 
\end{thm}

Theorem \ref{main} immediately implies the following corollary:

\begin{corollary}\label{cor}
If $t\geq 18$, $K_t^n$ is path-pairable. 
\end{corollary}

The above construction provides examples of path-pairable graphs on $N=t^n$ vertices with $\Delta = t\cdot n= \log N\cdot \frac{t}{\log t}$ maximum degree. Observe that $t$ can be chosen to be a constant ($t=18$) thus we have obtained path-pairable graphs on $N$ vertices with $\Delta\approx 4.3\log N$. 

Before the proof of Theorem \ref{main} we set the notation and terminology: Let us fix a coordinate (say, the last coordinate) in the $n$-dimensional grid $K_t^n$ and partition the vertices into $n$ copies of $K_t^{n-1}$ with respect to their last coordinates. Observe that every partition defines a subgraph isomorphic to an $(n-1)$-dimensional grid $K_t^{n-1}$. We denote these subgraphs by $L_1,\dots,L_{t}$ and call them {\it layers}. Similarly, fixing the first $n-1$ coordinates results in $t^{n-1}$ copies of $K_t$; we denote these complete subgraphs by $l_1,\dots, l_{t^{n-1}}$ and refer to them as {\it columns}.

\section{Proof of Theorem \ref{main}} We prove our statement by an inductive approach; given an edge $uv$ of the demand graph with $u,v\in K_t^n$ we replace $uv$ by a path of three edges $uu'$, $u'v'$, and $v'v$ where $u',v'\in K_t^n$ and $u,u'$ and $v,v'$ lie in the same columns and $u',v'$ share the same layer, that is, $u',v'\in L_i$ for some $i\in\{1,2,\dots,t\}$. Having done that we consider the new demand edges defined within the $t$ layers and $t^{n-1}$ columns and break the initial problem into $t^{n-1}+t$ subproblems that we solve inductively. We devote the upcoming sections to the detailed discussion of the above described solution plan. 

For the discussion of the base case $n=1$ as well as for the inductive step we use the following theorem:

\begin{thm}{\cite{tpc}}\label{delta_n/3}
Let $I(G)= K_t$ and add $q$ terminal vertices to every vertex of $I(G)$; we denote this graph by $K_t(q)$. If $q\leq 2\lfloor\frac{t}{6}\rfloor-2$, then $K_t(q)$ is terminal-pairable.
\end{thm}

We mention that instead of using Theorem \ref{delta_n/3} we could use a weaker version of the theorem with $q\leq \frac{t}{4+2\sqrt{3}}$ proved by Csaba et al. in \cite{CS}. With every further step of our proof unchanged a result similar to Theorem \ref{main} could be proved with a smaller bound on $\Delta(D)$ and a worse (larger) constant multiplier but same order of magnitude on $\Delta$ in our path-pairable graphs.

Let $q$ be an even number with $2\leq q \leq \lfloor\frac{t}{6}\rfloor-1$ and let $D=(V(D),E(D))$ be a demand multigraph with $V(D)=K_t^n$ and $\Delta(D)\leq q$. Let $E'(D)$ denote the set of demand edges whose endvertices lie in different $l_i, l_j$ columns. We construct an auxiliary graph $H$ with $V(H)=V(K_t^{n-1})$ and project every edge of $E'(D)$ into $H$ by deleting the last coordinates of the endvertices. It is easy to see that $\Delta(H)\leq t\cdot q $. We may assume without loss of generality that $D$ is $t\cdot q$-regular by joining additional pairs of vertices or replacing edges by paths of length two if necessary. We use the well known 2-Factor-Decomposition-Theorem of Petersen \cite{Petersen} to distribute the original demand edges among the layers $L_1,\dots, L_t$ and define new subproblems on them:

\begin{thm}[\cite{Petersen}]\label{Petersen}
Let $G$ be a $2k$-regular multigraph. Then $E(G)$ can be decomposed into the union of $k$ edge-disjoint $2$-factors.
\end{thm}

Obviously, the graph $H$ satisfies the conditions of Theorem \ref{Petersen} thus $E(H)$ can be partitioned into $\frac{q}{2}\cdot t$ edge-disjoint 2-factors. By arbitrarily grouping the above two factors into $\frac{q}{2}$-tuples we can partition $E(H)$ into $t$ edge disjoint subgraphs $H_1,\dots,H_t$ with $\Delta(H_i)\leq q$.

Assume now that the vertices $u=(\underline{a},i)$ and $v=(\underline{b},j)$  ($a,b\in [t]^{n-1}$) are joined by a demand edge belonging to $E'(D)$ (thus $\underline{a}\neq \underline{b}$) and assume that the corresponding edge in $H$ is contained by $H_k$. We then replace the demand edge $uv$ by the following triple of newly established demand edges:
$(\underline{a},i)(\underline{a},k)$, $(\underline{a},k)(\underline{b},k)$, and 
$(\underline{b},k)(\underline{b},j)$. We claim the following: 

\begin{itemize}
\item[i)] For  every layer $L_j$ the condition $\Delta(L_j)\leq q$ holds.
\item[ii)] For every layer $l_j$ the condition $\Delta(l_j)\leq 2q$ holds.
\end{itemize}
The first statement obviously follows from the partition of $E'(D)$. For the second one observe that a vertex $v$ in $l_j$ initially was incident to $q$ demand edges and at most $q$ additional demand edges have been joined to it (otherwise i) is violated). Observe now that every layer $L_j$ contains an $(n-1)$-dimensional subproblem that can be solved (within the layer) by the inductive hypothesis. Also, every layer $l_j$ contains a subproblem that can be solved by the Theorem \ref{delta_n/3}. If $n=2$ we also use the result of the above Theorem instead of the inductive hypothesis. That completes our proof.

\section{Conclusions and additional remarks}

By using Theorem \ref{delta_n/3} and the described inductive approach we proved that $K_t^n$ is path-pairable for $t\geq 18$, $n\in\mathbb{Z}^+$.
It was conjectured by Faudree, Gy\'arfas, and Lehel \cite{pp} that the result of Theorem \ref{delta_n/3} is true for $q\leq \lfloor\frac{t}{2}\rfloor$. If the conjecture is true it improves the constant 4.3 and decreases the lower bound on $t$ in Corollary \ref{cor}, yet it does not effect the order of magnitude $\log N$ of $\Delta$.

We mention that one particularly interesting and promising path-pairable candidate (with the same order of magnitude but better constant for $\Delta$) is the $n$-dimensional hypercube $Q_n$ on $2^n$ vertices ($\Delta(Q_n) = n$). Observe that hypercubes are special members of the above studied complete grid family as $Q_n =(K_2)^n$.
Although it is known that $Q_n$ is not path-pairable
for even values of $n$ \cite{F}, the question is open for odd dimensional
hypercubes if $n\geq 5$ ($Q_1$ and $Q_3$ are both path-pairable).

\begin{conjecture}[\cite{CS}]
The $(2k+1)$-dimensional hypercube $Q_{2k+1}$ is path-pairable for all
$k\in\mathbb{N}$.
\end{conjecture}
\section*{Acknowledgment}
We would like to heartfully thank Professor Jen\H o Lehel for drawing our attention to the above discussed problems. The third author in particular wishes to express his gratitude for the intriguing discussions about terminal-pairability problems.

\bibliographystyle{acm}
\bibliography{refs.bib}

\end{document}